\documentclass{amsart}
\usepackage{amsfonts}

\newtheorem{theorem}{Theorem}
\theoremstyle{plain}

\newtheorem{corollary}{Corollary}

\newtheorem{definition}{Definition}

\newtheorem{lemma}{Lemma}

\newtheorem{remark}{Remark}

\numberwithin{equation}{section}

\begin{document}
\title[A new family of polynomial identities for computing determinants]{A
new family of polynomial identities for computing determinants}
\author{Georgy\ P.\ Egorychev}
\address{RUSSIA,\ Krasnoyarsk,}
\email{anott@scn.ru}
\date{December 3, 2011}

\begin{abstract}
We give new definitions for the determinant over commutative ring
$K$, noncommutative ring $\mathbf{K}$, noncommutative ring
$\mathcal{K}$ with associative powers, over noncommutative nonassociative
ring $\mathfrak{K}$, and study their properties.
\end{abstract}

\maketitle

Let $K$ be a commutative ring, $\mathbf{K}$ a noncommutative associative
ring, $\mathcal{K}$ a noncommutative ring with associative powers
(one-monomial associativity), and $\mathfrak{K}$ be a noncommutative
nonassociative ring; let each ring be with division by integers. Here we
obtain a new family of polynomial identities for determinant over the ring
$K$, which allows us to give new definitions for determinants over rings
$\mathbf{K}$, $\mathcal{K}$, $\mathfrak{K}$, and to study their
properties. These definitions are closely related to the definition of
symmetrized \textit{Barvinok's} determinant $sdet(A)$ \cite{Bar2000}
different from the well-known determinant of Dieudonn\'{e} over a division
ring, the quasideterminant \cite{GelRet1991}, and other well-known
determinants over noncommutative associative rings \cite{Arv2009}. We also
estimate the computational complexity of the obtained formulas for the
determinants.

Let $A=\left( a_{ij}\right) $ be an $n\times n$ matrix with elements from
the ring $K$. Let $S_{n}$ be the set of all permutations $\sigma =\left(
\sigma \left( 1\right) ,\ldots ,\sigma \left( n\right) \right) $ of the set $%
\{1,\ldots ,n\},$ $\tau (\sigma )$ be the number of inversions in $\sigma .$
Let $S_{n}^{\left( e\right) }$ and $S_{n}^{\left( o\right) }$ be the subsets
of even and odd permutations in $S_{n}$, respectively$.$ We call the
sequence of elements of $a_{1\sigma \left( 1\right) },\ldots ,a_{n\sigma
\left( n\right) }$ the diagonal $l\left( \sigma \right) $ of matrix $A,$ the
sequence of elements $a_{i_{1}\sigma \left( i_{1}\right) },\ldots
,a_{i_{k}\sigma \left( i_{k}\right) },$ $1\leq i_{1}<\ldots <i_{k}\leq n$
the subdiagonal $l$ of length $k$ of matrix $A$. We denote by $L_{k}^{\left(
e\right) }\left( L_{k}^{\left( o\right) }\right) $ the set of all
subdiagonals of length $k$ of diagonals from $S_{n}^{\left( e\right) }\left(
S_{n}^{\left( o\right) }\right) ,$ $k=1,\ldots ,n$. The function $su\left(
l\right) $ on the subdiagonal (diagonal) $l=diag(a_{i_{1}\sigma \left(
i_{1}\right) },\ldots ,a_{i_{k}\sigma \left( i_{k}\right) }),$ is defined to
be the sum of its elements, i.e. $su\left( l\right)
:=\sum_{s=1}^{k}a_{i_{s}\sigma \left( s\right) }.$

\begin{definition}
\label{D1}(the determinant over a commutative ring)%
\begin{equation*}
det\left( A\right) :=\sum_{\sigma \in S_{n}}\left( -1\right) ^{\tau \left(
\sigma \right) }a_{1\sigma \left( 1\right) }\ldots a_{n\sigma \left(
n\right) }.
\end{equation*}
\end{definition}

In the following Main Theorem we employ the well-known polarization theorem
\cite{Car,Egor2007} to obtain a new family of polynomial identities for
determinants that contain up to $n!$ free variables.

\begin{theorem}[Main Theorem]
\textit{If }$A=\left( a _{ij}\right) $\textit{\ is an }$n\times n$%
\textit{\ matrix over the ring }$K$\textit{, then the following formulas
hold}:%
\begin{equation}
\det \left( A\right) =\frac{1}{n!}\sum_{\sigma \in S_{n}}\left( -1\right)
^{\tau \left( \sigma \right) }\{\left( -1\right) ^{n}\gamma _{\sigma
}^{n}+\sum_{k=1}^{n}(\left( -1\right) ^{n-k}\sum_{1\leq j_{1}<\ldots
,j_{k}\leq n}(\gamma _{\sigma }+\sum_{s=1}^{k}a_{_{j_{s}}\sigma \left(
_{j_{s}}\right) })^{n})\},  \label{B3}
\end{equation}%
\textit{where }$\{\gamma _{\sigma }\}_{\sigma \in S_{n}}$ \textit{is the set
of }$n!$ \textit{free variables }$\gamma _{\sigma }\in K${\Large $\mathbf{.}$
}\textit{In particular for} \textit{all} $\gamma _{\sigma }=\gamma $ \textit{%
and for all }$\gamma _{\sigma }=0$ \textit{we have, correspondingly}%
\begin{equation*}
det\left( A\right) =\frac{\left( -1\right) ^{n-1}}{n!}\{\sum_{l\in
L_{n-1}^{\left( e\right) }}(\gamma +su\left( l\right) )^{n}-\sum_{l\in
L_{n-1}^{\left( o\right) }}(\gamma +su\left( l\right) )^{n}+
\end{equation*}%
\begin{equation}
+\sum_{l\in L_{n}^{\left( o\right) }}(\gamma +su\left( l\right)
)^{n}-\sum_{l\in L_{n}^{\left( e\right) }}(\gamma +su\left( l\right) )^{n}\},
\label{B4}
\end{equation}%
\begin{equation}
det\left( A\right) =\frac{\left( -1\right) ^{n-1}}{n!}\{\sum_{l\in
L_{n-1}^{\left( e\right) }}su^{n}\left( l\right) -\sum_{l\in L_{n-1}^{\left(
o\right) }}su^{n}\left( l\right) +\sum_{l\in L_{n}^{\left( o\right)
}}su^{n}\left( l\right) -\sum_{l\in L_{n}^{\left( e\right) }}su^{n}\left(
l\right) \}.  \label{B5}
\end{equation}
\end{theorem}

\begin{corollary}
\label{C1}If $A=\left( a _{ij}\right) $\ \ is an $n\times n$\ matrix
over $\mathbb{C}
$\ then for $t=1,2,\ldots ,n-1$\ the following identities are hold:%
\begin{equation}
\sum_{l\in L_{n-1}^{\left( e\right) }}su^{t}\left( l\right) -\sum_{l\in
L_{n-1}^{\left( o\right) }}su^{t}\left( l\right) +\sum_{l\in L_{n}^{\left(
o\right) }}su^{t}\left( l\right) -\sum_{l\in L_{n}^{\left( e\right)
}}su^{t}\left( l\right) =0.  \label{B6}
\end{equation}
\end{corollary}

\begin{corollary}
\label{C2}(a new criterion for independence of matrix rows (columns)).%
\textit{\ If }$A=\left( a_{ij}\right) $\textit{\ is an }$n\times n$%
\textit{\ matrix over the ring }$K,$\textit{\ then }$det\left( A\right) =0$%
\textit{\ iff the following identity holds}%
\begin{equation*}
\sum_{l\in L_{n}^{\left( e\right) }}su^{n}\left( l\right) -\sum_{l\in
L_{n}^{\left( o\right) }}su^{n}\left( l\right) =\sum_{l\in L_{n-1}^{\left(
e\right) }}su^{n}\left( l\right) -\sum_{l\in L_{n-1}^{\left( o\right)
}}su^{n}\left( l\right) .
\end{equation*}
\end{corollary}

\begin{remark}
\textit{a})\textit{\ The formulas }$(\ref{B3})$ \textit{-- }$(\ref{B5})$%
\textit{\ are obtained using all characteristic properties of the
determinant }\cite{Muir}.\newline
\textit{b})\textit{\ Formula }$(\ref{B3})$\textit{\ for }$n!$\textit{\ free
variables }$\{\gamma _{\sigma }\}_{\sigma \in S_{n}}$ \textit{generates }$%
2^{n!}$\textit{\ different polynomial identities if we put each }$\gamma
_{\sigma }=0$ or $\gamma _{\sigma }\neq 0$\textit{. Each of these formulas
can be taken as the definition for }$det\left( A\right) ,$\textit{\ it
requires its own number of arithmetical operations when calculated and
according to Corollary \ref{C1} generates a new set of identities.}\newline
\textit{c}) \textit{Formula} $(\ref{B3})$\textit{\ and its special cases
employ (besides divisions by $n!$)} \textit{\ only operations }$+$\textit{,}$%
-$\textit{,\ and raising to power of }$n$\textit{, yet does not use a
commutativity of multiplication in the ring} $K$.
\end{remark}

\begin{definition}
\label{D2}(the noncommutative determinant $edet(A)$ over ring $\mathbf{K}$
or $\mathcal{K}).$%
\begin{equation}
edet(A):=\frac{\left( -1\right) ^{n-1}}{n!}\{\sum_{l\in L_{n-1}^{\left(
e\right) }}su^{n}\left( l\right) -\sum_{l\in L_{n-1}^{\left( o\right)
}}su^{n}\left( l\right) +\sum_{l\in L_{n}^{\left( o\right) }}su^{n}\left(
l\right) -\sum_{l\in L_{n}^{\left( e\right) }}su^{n}\left( l\right) \}.
\label{B9a}
\end{equation}
\end{definition}

\begin{lemma}
\label{L3}(the properties  of$edet(A)$ over ring $\mathbf{K}$).\newline
$a)$ $edet(A)$ \textit{is a polyadditive} \textit{function} \textit{of the rows
and columns of matrix} $A;$\newline
$b)$ $edet(A)\mathit{\ }$\textit{is an antisymmetric function of the rows and
columns of matrix }$A$;\newline
$c)$ $edet(A)=0,$ \textit{\ if the} \textit{matrix }$A\ $\textit{contains
zero row or zero column};\newline
$d)$ $edet(A)=0,$ \textit{if the} \textit{matrix }$A\ $\textit{contains two
equal rows }(\textit{columns});\newline
$e)$ $edet(A)=edet(A^{T}).$\newline
\textit{The Laplace formulas and the formula }$\det (AB)$ = $\det (A)\times
\det (B)$ \textit{are not valid in the general case.}
\end{lemma}

Let $sdet\left( A\right) $\textit{\ }be the symmetrized Barvinok's
determinant over ring\textit{\ }$\mathbf{K}$ \cite{Bar2000}:%

\begin{equation*}
sdet(A):=\frac{1}{n!}\sum_{\mu \in S_{n}}\sum_{\sigma \in S_{n}}\left(
-1\right) ^{\tau \left( \sigma \right) +\tau\left( \mu \right)}a_{\mu \left( 1\right) \sigma \left(
1\right) }\ldots a_{\mu \left( n\right) \sigma \left( n\right) }.
\end{equation*}

\begin{lemma}
\label{L4}\textit{The following formulas are valid for }$edet\left( A\right)
$ \textit{over ring} $\mathbf{K}$:%
\begin{equation*}
sdet(A)=edet(A),
\end{equation*}%
\begin{equation*}
edet(I_{n})=e,\text{\textit{\ if the ring} }\mathbf{K}\text{ \textit{%
contains the unit} }e\text{, \textit{and} }I_{n}\text{ \textit{is the
identity} }n\times n\text{ \textit{matrix}}.
\end{equation*}
\end{lemma}

\begin{lemma}
\label{L5}(computational complexity of $edet(A)$). \textit{The
computation of} $edet(A)$\textit{\ via formulas }$(\ref{B5})$\textit{\ and} $(%
\ref{B9a})$\textit{\ for large }$n$ \textit{requires of order\ }$%
(3n-1)\times n!$\ \textit{additions and }$(n+1)!\times \ln n$\ \textit{%
mulltiplications} (see \cite{Koch1994})
\end{lemma}

Formula $(\ref{B9a})$ allows to introduce the following new general
definition

\begin{definition}
\label{D3}(the determinant $edet(A)$ over ring $\mathfrak{K}$).%
\begin{equation*}
edet(A):=\frac{\left( -1\right) ^{n-1}}{n!}\{\sum_{l\in L_{n-1}^{\left(
e\right) }}Ass(su^{n}\left( l\right) )-\sum_{l\in L_{n-1}^{\left( o\right)
}}Ass(su^{n}\left( l\right) )+\sum_{l\in L_{n}^{\left( o\right)
}}Ass(su^{n}\left( l\right) )-\sum_{l\in L_{n}^{\left( e\right)
}}Ass(su^{n}\left( l\right) )\},
\end{equation*}%
\textit{where }$Ass(a^{n})$, $a\in $\textit{\ }$\mathfrak{K}$, \textit{is}
\textit{associate operator}%
\begin{equation*}
Ass(a^{n}):=\{(a\times (a\times (a\ldots (a\times a))))+\ldots +((((a\times
a)\times a)\times a)\ldots \times a)\}/C_{n},
\end{equation*}%
\textit{and }$C_{n}=\left( 2n-2\right) !/n!\left( n-1\right) !$\textit{\ is
the Catalan number of }$n$\textit{-products in nonassociative algebra }(%
\textit{all arrangements of brackets in the }$n$-\textit{product} $a\times
a \times \ldots \times a )$\textit{.}
\end{definition}

\begin{lemma}
\label{L6}(the properties of the determinant $edet(A)$ over nonassociative
rings $\mathcal{K}$ and $\mathfrak{K}$). \textit{The determinant }$edet(A)$%
\textit{\ over rings }$\mathcal{K}$ \textit{and} $\mathfrak{K}$\textit{\
satisfies properties b) -- e)\ of Lemma \ref{L3}. The Laplace formulas, the
polyadditivity by rows and columns of matrix} $A$ \textit{and the formula} $%
edet(AB)=edet(A)\times edet(B)$ are not valid in general case.\textit{\ }
\end{lemma}

Obtaining similar results is of interest for Shur functions, the mixed
discriminants and many other matrix functions of planar and space matrices
related to determinants. Our results may find applications in the theory of
permanents \cite{Egor2007}, the theory of $n$-Lie algebras \cite{Fil},
differential geometry {\cite{Sab}} and elsewhere.

The author is grateful to his colleagues M.N. Davletshin, V.M. Kopytov, Y.N.
Nuzhin, A.V. Timofeenko, I.P. Shestakov for discussion of the basic results
of this work and a number of useful remarks.


\begin{thebibliography}{99}
\bibitem{Arv2009} Arvind V. and Srinivasa S. On the Hardness of
Noncommutative Daterminant, Electronic Colloquium on Computational
Camplexity, Report No. 103, 2009.

\bibitem{Bar2000} Barvinok A. New permanent estimators via non-commutative determinants. Arxiv preprint math/0007153, 2000, arXiv:math/0007153.

\bibitem{Car} Cartan H. Elementary theory of analytic functions of one or
several complex variables. Dover Publ., New York, 1995, 228 p.

\bibitem{Egor2007} Egorychev G.P. Discrete mathematics. Permanents.
Krasnoyarsk, Siberian Federal University, 2007, 272 p., (transl. in English:
Springer, 2011).

\bibitem{Fil} Filippov V.T. On the n-Lie algebra of Jacobians, Sibirsk. Mat. Zh., 39:3, 1998, 660--669.

\bibitem{GelRet1991} Gel'fand I. M. and Retakh V. S. Determinants of matrices over noncommutative rings, Funktsional. Anal. i Prilozhen. 25, 1991, No. 2, 13--25, 96.

\bibitem{GevSlinShes82} Zevlakov K.A., Slin'ko A.M., Shestakov I.P. and
Shirshov A.I. Rings that are nearly associative, Acad. Press, 1982, 432 p.

\bibitem{Koch1994} Kochergin V.V. About complexity of computation one-terms of powers, Discrete analysis, IM SO RAN, v. 27, Novosibirsk, 1994, 94 -- 107 (in Russian)

\bibitem{Muir} Muir T. A Treatise on the Theory of Determinants, New York:
Dover, 1960.

\bibitem{Sab} Sabinin L.V. Methods of Nonassociative Algebra in Differential
Geometry. in Supplement to Russian translation of S.K. Kobayashi and K.
Nomizu "Foundations of Differential Geometry", 1, Moscow, Nauka, 1981.
\end{thebibliography}
\end{document}